\documentclass[11pt]{article}
\usepackage{amsfonts}
\usepackage{mathrsfs}
\usepackage{graphics}
\usepackage{indentfirst}
\usepackage{cite}
\usepackage{latexsym}
\usepackage{amsmath}
\usepackage{amssymb}
\usepackage[dvips]{epsfig}
\usepackage{amscd}
\newcommand{\ba}{\begin{aligned}}
\newcommand{\ea}{\end{aligned}}
\newcommand{\be}{\begin{equation}}
\newcommand{\ee}{\end{equation}}
\newcommand{\bnn}{\begin{eqnarray*}}
\newcommand{\enn}{\end{eqnarray*}}

\renewcommand{\div}{ {\rm div }  }

\newcommand{\na}{\nabla }

\newcommand{\pa}{\partial}

\newcommand{\bl}{\begin{lemma}}
\newcommand{\el}{\end{lemma}}
\newcommand{\et}{\end{theorem}}

\newcommand{\te}{\theta}

\newcommand{\de}{\delta}

\newcommand{\la}{\label}
\newcommand{\ka}{\kappa}

\newcommand{\n}{\rho}

\def\la{\label}
\def\na{\nabla}

 \def\p{\partial}

\def\norm[#1]#2{\|#2\|_{#1}}
\numberwithin{equation}{section}

 \usepackage{indentfirst}
\usepackage{graphics}
\usepackage{epsfig}

\newtheorem{claim}{\bf \t}[part]

\newtheorem{theorem}{Theorem}[section]
\newtheorem{corollary}[theorem]{Corollary}
\newtheorem{lemma}[theorem]{Lemma}
\newtheorem{proposition}[theorem]{Proposition}
\newtheorem{remark}{Remark}[section]

\def\t{\theta}
\def\k{\kappa}

\def\m{\mu}

\def\l{\lambda}

\def\r{\rho}

\def\f{\frac}
\begin{document}

\title{One New Blowup Criterion for the 2D Full  Compressible Navier-Stokes System}

\author{ Yun Wang
\thanks{School of Mathematics, Soochow University, 1 Shizi Street, Suzhou, Jiangsu 215006, China ({\tt ywang3@suda.edu.cn}).}
 }
\date{ }
\maketitle

\begin{abstract}

We establish a blowup criterion for the
two-dimensional (2D)  full  compressible Navier-Stokes system. The criterion is given in terms 
of the divergence of the velocity field only, and is independent of the temperature. The criterion tells that 
once the strong solution blows up, the $L^\infty$-norm for the divergence of velocity blows up.

\


\

AMS: 35Q35, 35B65, 76N10
\end{abstract}

\section{Introduction}

The motion of compressible viscous, heat-conductive, ideal polytropic fluid is governed by the full compressible Navier-Stokes system. Suppose that the 
domain occupied by the fluid is $\Omega$. The whole system on $(0, T)\times \Omega$ consists of the following equations
\begin{eqnarray} \label{1.1}
\begin{cases}
\r_{t}+\mbox{div}(\r u)=0,\\
(\r u)_{t}+\mbox{div}(\r u\otimes
u)-\mu\Delta{u}-(\mu+\l)\nabla\mbox{div}u
+\nabla{P}=0,\\
c_v[(\r \t)_{t}+\mbox{div}(\r u\t
)]-\k\Delta{\t}+P\mbox{div}u=2\mu|\mathfrak{D}(u)|^2+\l(\mbox{div}u)^2.
\end{cases}
\end{eqnarray}
together with the initial-boundary conditions
\begin{eqnarray}
&(\rho, u, \theta)|_{t=0} = (\rho_0, u_0, \theta_0), \label{initial}\\
&u=0, \ \ \ \ \frac{\partial \theta}{\partial \vec{n}} = 0,\ \ \ \ \mbox{on}\ \partial \Omega. \label{boundary}
\end{eqnarray}
In this paper, we consider the 
two-dimensional case, i.e.,  $\Omega$ is a bounded smooth domain in $\mathbb{R}^2$. 
 $\n, u=\left(u^1,u^2\right)^{\rm tr},  $ $\te $ and $P=R\r\t \,(R>0) $
represent respectively the fluid density, velocity,  absolute
temperature and pressure. In addition, $ \mathfrak{D}(u)$ is the
deformation tensor: \bnn \la{z1.2}   \mathfrak{D}(u)  =
\frac{1}{2}(\nabla u + (\nabla u)^{\rm tr}). \enn
 The constant
viscosity coefficients $\mu$ and $\lambda$  satisfy the physical
restrictions:
 \be\la{h3} \mu>0,\quad \mu + \lambda\ge 0.
\ee Positive constants $c_v$  and  $\ka$ are respectively  the heat capacity and the ratio of the heat conductivity coefficient over the heat capacity.

There is a huge amout of literature on the existence and large time behavior of solutions to \eqref{1.1}. For the case that the initial density is far away from vacuum, local existence and uniqueness 
of classical solutions were proved in \cite{Nash, Serrin-1, Kazhikhov}. Matsumura-Nishida\cite{MN} first obtained the global classical solution when the initial data is close to a non-vacuum equilibrium in some 
Sobolev space $H^s$. Later, Hoff\cite{Hoff-JDE} constructed a global weak solution for discontinuous initial data, with the assumption that the initial density is strictly positive. 

Normally  the presence of vacuum state makes the problem more complicated. We recall some results about the weak solution for this case first. For the global weak solution to the barotropic case, the major breakthrough is due to Lions\cite{Lions} and subsequently improved by Jiang-Zhang\cite{JZ} and Feireisl\cite{Feireisl}. 
They succeeded in constructing a global weak solution with finite energy when the pressure $P = \rho^\gamma, \ \gamma>\frac{N}{2}$, where $N$ is the dimension. Recently, Huang-Li\cite{HL-full} obtained the global weak solution to the full Navier-Stokes system \eqref{1.1} provided the initial energy is suitably small. Strong solutions have also been under investigation. The first local existence result was derived by Cho-Kim\cite{Cho-Kim}. For a global classical solution with small energy, refer to \cite{HLX-CPAM, HL-full}. On the other hand, Xin\cite{Xin} first contributed a remarkable blow-up result. He proved that if the initial density has compact support, any smooth solution to the Cauchy problem of the full Navier-Stokes equations without heat conduction blows up in finite time. In this direction, for more recent progress, see \cite{Cho-Jin, Rozanova} and references therein. 

Taking into consideration both the local existence results and the blowup results, then it is important to study the mechanism of possible blowup and structure of possible singularity. For the blowup criterion for the  compressible flow, there have been many results, \cite{Cho-Jin, FJO, FZZ, Huang, HL-MAA, HLW, HLX-SIAM, HLX-CMP, HX, JO, SWZ-JMPA, SWZ-ARMA, Wen-Zhu}. It should been mentioned here that Huang-Li-Xin\cite{HLX-CMP} first established a Beale-Kato-Majda type blowup criterion for the baratropic case. In fact they proved that if $T^*$ is the maximal time of existence for local strong solution, then 
\begin{equation}\label{HLX}
\lim_{T\rightarrow T^*} \int_0^T \|\nabla u\|_{L^\infty} \, dt = \infty,
\ee
under the assumption $7\mu > \lambda$ when $\Omega$ is a three-dimensional domain. Jiang-Ou\cite{JO} extended this criterion to the full Navier-Stokes system \eqref{1.1} over a periodic domain or unit square domain of $\mathbb{R}^2$ and proved that 
\begin{equation}\label{JO}
\lim_{T\rightarrow T^*} \int_0^T  \|\nabla u\|_{L^\infty} \, dt = \infty. 
\ee
Recently, Huang-Li-Wang\cite{HLW} obtained a Serrin type blow up criterion for \eqref{1.1} in $\mathbb{R}^N$. Here is the criterion, 
\be\la{HLW}
\lim_{T\rightarrow T^*} \int_0^T  \left( \|{\rm div} u\|_{L^\infty}  + \|u\|_{L^r}^s \, dt  \right) = \infty,\ \ \ \ \frac{2}{s}+ \frac{N}{r} =1, \ \ N<r \leq \infty.
\ee
which is analogue to the Serrin criterion for the 3D incompressible Navier-Stokes equations. In particular, for $N=2$, if one can bound a priorily $\|u\|_{L^2(0, T; L^\infty)}$-norm or $\|u\|_{L^4(0, T; L^4)}$-norm, then 
\eqref{HLW} can be replaced by 
\be\label{divergence}
\lim_{T \rightarrow T^*} \int_0^T \|{\rm div} u\|_{L^\infty} dt = \infty. 
\ee
If \eqref{divergence} is proved, it is an improvement of the work by Jiang-Ou\cite{JO} and it reveals some connection between the compressible and incompressible Navier-Stokes equations, since global strong solutionwith vacuum has been proved for 2D incompressible case\cite{HW}. The question is we can not get the uniform bound of $\|u\|_{L^2(0, T; L^\infty)}$ or $\|u\|_{L^4(0, T; L^4)}$ from the a priori energy estimate. The aim of our paper is to verify 
\eqref{divergence} and the key trick is the use of Lemma \ref{logarithmic-Sobolev} below, one critical Sobolev embedding inequality. 

The results such as \eqref{HLX} or \eqref{JO} or \eqref{HLW}, notice that they are all in terms of velocity field only. There is another big class of results which are in terms of density $\rho$ and temprature $\theta$. For example, 
Sun-Wang-Zhang\cite{SWZ-JMPA} obtained the following criterion in 3D, 
\be\label{SWZ}
\lim_{T\rightarrow T^*} \left(  \sup_{0\leq t \leq T}  \{ \|\rho\|_{L^\infty} + \|\rho^{-1}\|_{L^\infty} + \|\theta\|_{L^\infty} \} \right) = \infty.
\ee
Fang-Zi-Zhang\cite{FZZ} extended the result to the 2D problem with a refiner form, 
\be\label{FZZ}
\lim_{T\rightarrow T^*} \sup_{0\leq t \leq T} \{\|\rho\|_{L^\infty} + \|\theta\|_{L^\infty} \} = \infty. 
\ee

 Before stating our main result, we first explain the
notations and conventions used throughout this paper. We denote
$$\int fdx=\int_{\Omega}fdx.$$
For $1\le p\le \infty$ and integer $k\ge 0$, the standard homogeneous and inhomogeneous Sobolev
spaces are denoted by:
   \bnn \begin{cases} L^p=L^p(\Omega),\quad W^{k,p}=W^{k,p}(\Omega),\quad
   H^k=W^{k,2},\\
W_0^{1, p} =\{u \in W^{1, p} | \ u= 0\ \mbox{on}\ \partial \Omega\}, \ \ H_0^1 = W_0^{1, 2}.
 \end{cases}\enn
Let $$\dot{f} := f_t + u \cdot \nabla f$$
denote the material derivative of $f$.

Since we are going to work with the blowup criterion of the strong solutions, we'd like to recall the 
result for the existence of the local strong solution. The solution to the 3D full Navier-Stokes system with vacuum 
was obtained by Cho-Kim\cite{Cho-Kim}. The method there can be applied to the case in this paper, i.e. the case that $\Omega$ is a 
bounded domain in $\mathbb{R}^2$. And the corresponding result can be stated as follows(or refer to \cite{FZZ}):

\begin{theorem}\label{local-existence}
Let $q \in (2, \infty) $ be a fixed constant. Assume that the initial data satisfy 
$$\rho_0 \geq 0, \ \ \ \rho_0 \in W^{1, q}, \ \ u_0 \in H_0^1\cap H^2, \ \ \theta_0 \in H^2,$$
with the compatibility conditions
\begin{eqnarray} 
\mu \Delta u_0 + (\lambda + \mu ) \nabla {\rm div} u_0 - R \nabla (\rho_0 \theta_0) = \rho_0^{\frac12} g_1, \label{compatible-1}\\
\kappa \Delta \theta_0 + \frac{\mu}{2} |\nabla u_0 + (\nabla u_0)^{tr}|^2 + \lambda ({\rm div}u_0)^2 = \rho_0^{\frac12}g_2, \label{compatible-2}
\end{eqnarray}
for some $g_1, g_2\in L^2$. Then there exist a positive constant $T_0$ and a unique strong solution $(\rho, u, \theta)$ to the system
(\ref{1.1})-(\ref{boundary})such that 
\begin{eqnarray}
&\rho \geq 0,\ \ \ \rho \in C([0, T_0]; W^{1, q}),\ \ \ \theta\in C([0, T_0]; H^2), \\
&u \in C([0, T_0]; H_0^1 \cap H^2 ),\ \ \ \ (u, \theta) \in L^2 ([0, T_0]; W^{2, q}),\\
&(u_t, \theta_t) \in L^2 ([0, T_0]; H^1),\ \ (\sqrt{\rho}u_t, \sqrt{\rho}\theta_t ) \in L^\infty ([0, T_0]; L^2).
\end{eqnarray}
\end{theorem}

Regarding the blowup criterion for the local strong solution, here is our main theorem.

\begin{theorem}\label{thm1.1} Suppose the assumptions in  Theorem \ref{local-existence} are satisfied and $(\rho, u, \theta)$ is 
the strong solution. Let $T^*$ be the
maximal time of existence for that strong solution. If $T^* < \infty$, then
\begin{equation}\label{criterion}
\lim_{T\rightarrow T^*}
\|{\rm div}u\|_{L^1(0,T;L^{\infty})} = \infty .
 \end{equation}
   \end{theorem}

A few remarks are in order, 
\begin{remark}
It is worth noting that the conclusion in  Theorem \ref{thm1.1} is somewhat surprising  since the criterion \eqref{criterion} is independent of the temperature  and is the same as that of barotropic case(\cite{HLX-SIAM}).
In fact,  it seems that the nonlinearity of the highly nonlinear terms $|\mathfrak{D}(u)|^2$ and $(\div u)^2$ in the temperature equation
  is stronger than that of $\div (\n u\otimes u)$ in the momentum equations (\cite{SWZ-ARMA}), however,  \eqref{criterion}  shows that the nonlinear term $|\na u|^2$ can be controlled provided one can control  $\div (\n u\otimes u).$
\end{remark}

\begin{remark}
It is well known that the 2D incompressible homogenenous Navier-Stokes system has a unique global strong solution if the initial velocity belongs to $L^2$ or some more regular space, and recently it is proved in 
\cite{HW} that the 2D incompressible nonhomogenous Navier-Stokes system also has a unique global strong solution under some compatibility conditions, so the result in our paper is reasonable from this point. The blowup criterion here shows that ${\rm div}u$ plays an important role in the fluid dynamics. 
\end{remark}

\begin{remark}
The techniques in this paper can be easily adapted to the two dimensional periodic case. And the same criterion will be derived.
\end{remark}

The rest of the paper is organized as follows: In Section 2, we
collect some elementary facts and inequalities. The main result, Theorem \ref{thm1.1},  will be proved in Section 3.


\section{Preliminaries}
In this section, we recall some known facts and elementary
inequalities that will be used later. 

The first proposition is for the Lam\'e system, which comes from the momentum equation $(\ref{1.1})_2$. Assume that $\Omega\subset \mathbb{R}^2$ is a bounded 
smooth domain. Suppose $U\in H_0^1$ is a  weak solution to the Lam\'e system, 
\begin{eqnarray}\label{lame}
\begin{cases}
\mu \Delta U + (\mu + \lambda) \nabla {\rm div} U = F, \ \ \mbox{in}\ \Omega,\\
U(x) = 0,\ \ \ \mbox{on}\ \partial \Omega.
\end{cases}
\end{eqnarray}
In what follows, we denote $\mathcal{L}U= \mu \Delta U + (\mu + \lambda) \nabla {\rm div} U$. Owing to the uniqueness of solution, we denote $U = \mathcal{L}^{-1} F$. 

The system is an elliptic system under the assumption (\ref{h3}), hence some regularity estimates can be derived. For a proof, refer to \cite{SWZ-JMPA}.

\begin{proposition}\label{prop2.1} Let $q\in (1, \infty)$. Then there exists some constant $C$ depending only on $\lambda, \mu, p$ and $\Omega$ such that
\begin{itemize}
\item if $F\in L^p$, then
\begin{equation}\label{lame-2q}
\|U \|_{W^{2,p}} \leq C\|F\|_{L^p};
\end{equation}

\item if $F\in W^{-1, p}$(i.e., $F = {\rm div}f$ with $f = (f_{ij})_{2\times 2}$, $f_{ij}\in L^p $), then
\begin{equation}
\|U\|_{W^{1,p}} \leq C\|f\|_{L^p}.
\end{equation}

Moreover, for the endpoint case, if $f_{ij}\in L^\infty \cap L^2$, then $\nabla U \in BMO(\Omega)$ and there exists some constant C depending only on 
$\lambda, \mu, \Omega$ such that
\begin{equation}
\|\nabla U\|_{BMO(\Omega)} \leq C \left(\|f\|_{L^\infty} + \|f\|_{L^2}\right).
\end{equation}
Here $\|g\|_{BMO(\Omega)} := \|g\|_{L^2} + [g]_{BMO(\Omega)}, $ with
\begin{equation}\nonumber 
[g]_{BMO(\Omega)} : = \sup_{x\in \Omega,\ r\in (0, d)} \frac{1}{|\Omega_r(x)|} \int_{\Omega_r(x)} |g(y)- g_{\Omega_r(x)} | \, d y,
\end{equation}
\begin{equation}\nonumber
g_{\Omega_r(x)} = \frac{1}{|\Omega_r(x)|} \int_{\Omega_r(x)} g(y)\,  dy,
\end{equation}
where $\Omega_r(x) = B_r(x) \cap \Omega$ and $|\Omega_r(x)|$ denotes the Lebesgue measure of $\Omega_r(x)$. 
\end{itemize}

\end{proposition}


Two logarithmic Sobolev inequalities will be presented, which originate from the work owing to 
Brezis-Gallouet\cite{BG} and Brezis-Wainger\cite{BW}. The first one, together with Proposition \ref{prop2.1}, will give the estimate of $\nabla \rho$. For its proof, see also \cite{SWZ-JMPA}.
\begin{lemma}\label{lem2.2} Let $\Omega$ be a bounded Lipschitz domain in $\mathbb{R}^2$ and $f\in W^{1, p}$ with $p\in (2, \infty)$, there exists a constant $C$ depending only on $p$ such that
\begin{equation}
\|f\|_{L^\infty} \leq C \left( 1+  \|f\|_{BMO(\Omega)} \ln (e + \|f\|_{W^{1,p}})\right).
\end{equation}
\end{lemma}

The second inequality is in terms of both 
space integral and time integral. The proof can be found in \cite{HW} or refer to \cite{KOT} for a similar proof. It plays an important role in the proof of Lemma \ref{lem3.3}.
\begin{lemma}\label{logarithmic-Sobolev}
Let $\Omega$ be a smooth domain in $\mathbb{R}^2$, and $f\in L^2(0, T; H_0^1\cap W^{1, p})$, with $p>2$. Then there exists a constant $C$ depending only on $p$ such that
\begin{equation}\label{logarithmic}
\|f\|_{L^2(0, T; L^\infty)}^2 \leq C \left(1+ \|f\|_{L^2(0, T; H^1)}^2 \ln (e + \| f\|_{L^2(0, T; W^{1, p})}) \right) .
\end{equation}
\end{lemma}


\section{Proof  of Theorem \ref{thm1.1}.}

  Let $(\rho,u,\t)$ be a strong
solution described in Theorem \ref{thm1.1}. Suppose that
\eqref{criterion} is false, i.e., 
\begin{equation}\label{criterion-1}
\lim\limits_{T\rightarrow T^\ast}
\|{\rm div} u\|_{L^1(0,T;L^\infty)} 
\leq M_0<+\infty,
\end{equation}
which together with  $\eqref{1.1}_1$ yields immediately the following
 upper bound
of the density(see \cite[Lemma 3.4]{HLX-SIAM}).
\begin{lemma}\label{lem3.1}
 Assume that $\eqref{criterion-1}$ holds.
Then it is true that for $0\le T<T^*,$
 \begin{equation}\label{3-1-1}
\sup_{0\leq t\leq T}\|\r\|_{L^\infty}\leq C,
\end{equation}
 where and in what follows,
$C,$ $C_1$, $C_2$ and $C_3$  denote  generic constants depending only on $
  M_0,$ $ \mu,$ $\lambda, $  $R,\kappa,$ $c_v,$  $T^*,$  and the
initial data.
\end{lemma}


The next estimate is similar to an energy estimate. 
\begin{lemma}\label{lem3.2}
Under the assumption \eqref{criterion-1}, it holds that for $0\leq
T<T^\ast$,
\begin{eqnarray}\label{energy-estimate}
\sup_{0\leq t\leq T}\int\left(\r\t+\r |u|^2\right)  dx+\int_{0}^{T}\|\nabla
u\|_{L^2}^2\, dt\leq C.
\end{eqnarray}
\end{lemma}

  {\it Proof.} Applying standard maximum principle to $\eqref{1.1}_3 $ together with
   $\theta_0\ge 0$  (c.f. \cite{FJO,Feireisl}) shows that 
\begin{equation}\label{3-2-1} 
\inf_{\mathbb{R}^3\times [0,T]}\theta(x,t)\ge 0.\end{equation} 
 It follows from \eqref{1.1} that the specific energy $E\triangleq c_v\t+\f12|u|^2$ satisfies
\begin{equation}\label{3-2-2}
(\r E)_t+\mbox{div}(\r Eu+Pu)=\Delta(\k \t+\f12\m
|u|^2)+\m\mbox{div}(u\cdot\nabla u)+\l\mbox{div}(u\mbox{div}u).
\end{equation}
Integrating \eqref{3-2-2} over $\Omega \times [0,T]$ directly
yields
\begin{equation}\label{3-2-3}
\sup_{0\leq t\leq T}\int\left(\r\t+\r|u|^2\right)  dx\leq C.
\end{equation}
Multiplying $(\ref{1.1})_2$ by $u$  and integrating the resulting
equation over $\Omega$,  we obtain after using \eqref{3-2-1} and \eqref{3-2-3} that
\bnn\ba
&\f12\frac{d}{dt}\int\r|u|^2dx+\m\|\nabla u\|_{L^2}^2+(\m+\lambda)\|{\rm div} u\|_{L^2}^2 \\  & \ \ \leq  C
\|\mbox{div}u\|_{L^\infty}\int\r\t \, dx \\
&\ \ \leq 
C\|\mbox{div}u\|_{L^\infty},\ea
\enn
 which together with \eqref{criterion-1} and \eqref{3-2-3} gives \eqref{energy-estimate}.  It completes the proof  of Lemma \ref{lem3.2}.


\vspace{3mm}For a slightly higher order estimate, we derive the bound for the $L^\infty(0,T;
L^2)$-norm of $\nabla u,$  which will play a key role in obtaining more high order estimates. 
\begin{lemma}\label{lem3.3}
Under the condition \eqref{criterion-1}, it holds that for $0\leq
T<T^\ast$,
\begin{equation}\label{3-3-1}\sup_{0\leq t\leq T}
\int\left(\r\t^2 +|\nabla u|^2\right)
\, dx+\int_{0}^{T}\int\left(|\nabla\t|^2+\t|\nabla
u|^2+\r|\dot{u}|^2\right) dxdt \leq C.
\end{equation}
\end{lemma}


Before the proof of Lemma \ref{lem3.3}, we present an auxiliary lemma, which controls $L^p$-norm of $\theta$ by $\|\nabla \theta\|_{L^2}$. And it will help the proof of Lemma \ref{lem3.3}. 
\begin{lemma}\la{lem3.4}
Under the condition \eqref{criterion-1}, it holds  on $[0, T^*)$ that
for every $p\in [1, \infty)$,
\be \la{3-4-2}
\|\theta\|_{L^p} \leq C + C(p)\|\nabla \theta\|_{L^2}. 
\ee
\end{lemma}

{\it Proof.}\ Denote by $\bar{\theta} =\frac{1}{|\Omega|} \int \theta \, dx$, the average of $\theta$,  
\be\la{3-4-4} \ba 
|\bar{\theta}|  \int \rho \, dx  
& \leq \left|\int \rho \theta dx \right| + \left| \int \rho \left(\theta -\bar{\theta}\right) dx \right| \\
& \leq C + C \|\nabla \theta \|_{L^2},
\ea \ee
which together with Poincar\'e's inequality implies 
\be \la{3-4-5}
\|\theta\|_{L^2} \leq C + C\|\nabla \theta\|_{L^2}.
\ee
And consequently, \eqref{3-4-2} holds.

\vspace{3mm} {\it Proof of Lemma \ref{lem3.3}.} First, multiplying $(\ref{1.1})_3$ by $\t$ and
integrating the resulting equation over $\Omega$ lead to
\begin{equation}\label{3-3-2}
 \f{d}{dt}\int\r\t^2 \, dx+2\k\|\nabla\t\|_{L^2}^2\leq
C\|\mbox{div}u\|_{L^\infty}\int \r\t^2\, dx+C \int\t|\nabla u|^2\, dx.
\end{equation}

To make the estimate close, one needs to bound the term $\int \theta |\nabla u|^2 dx$ in \eqref{3-3-2}. To achieve that, we borrowed the idea from \cite{HLW}, 
 multiplying $\eqref{1.1}_2$ by $u\t$ and integrating the resulting
equation over $\Omega$. Then
\be\ba \label{3-3-3}
&\mu \int |\nabla u|^2 \theta\, dx + (\mu + \lambda) \int |{\rm div}u|^2 \theta \, dx  \\ = &
-\int \rho \dot{u} \cdot u\theta \, dx  - \mu \int u\cdot \nabla u \cdot \nabla \theta \, dx 
-(\mu + \lambda) \int {\rm div}u (u\cdot \nabla \theta) \, dx - \int \nabla P \cdot u \theta\,  dx.
\ea\ee

We estimate the terms on the righthand of \eqref{3-3-3}. By the Young's inequality,
\begin{equation}\label{3-3-4}
\left| \int \rho \dot{u}\cdot u \theta \, dx \right|\leq \eta \int \rho |\dot{u}|^2 \, dx + C_{\eta} \int \rho \theta^2 |u|^2\,  dx,
\end{equation}
and 
\be\ba \label{3-3-5}
& \left| \mu \int u\cdot \nabla u \cdot \nabla \theta \, dx 
+ (\mu + \lambda) \int {\rm div}u \cdot (u\cdot \nabla \theta) \, dx  \right| \\
& \ \ \ \leq \frac{\epsilon}{4} \|\nabla \theta\|_{L^2}^2 + C_{\epsilon} \int |u|^2 |\nabla u|^2\,  dx,
\ea \ee
where $\eta, \epsilon$ are  small positive constants to be determined later. Using integration by parts, 
\be \ba \label{3-3-6} 
\left|\int\nabla{P}u\t \, dx\right|&=\left|\int P\t\mbox{div}u\, dx+\int Pu\cdot\nabla\t
\, dx\right| \\
& \le \frac{\epsilon}{4} \|\nabla\t\|_{L^2}^2+C\|\mbox{div}u\|_{L^\infty}\int\r\t^2 \, dx
+C_{\epsilon} \int\r^2 \t^2|u|^2 \, dx.
\ea \ee

Combining the estimates \eqref{3-3-2}-\eqref{3-3-6}, we
obtain after choosing $\epsilon$ suitably small  that
\begin{eqnarray}\label{3-3-7}
&&\f{d}{dt}\int \r\t^2 \, dx+\int\left(\t|\nabla u|^2+ \k |\nabla\t|^2\right) \, dx\nonumber\\
\leq && C\eta\int\r|\dot u|^2  \, dx+
C\|{\rm div} u\|_{L^\infty}\int\r\t^2\, dx+ C_\eta  \int\left(\r\t^2 |u|^2 + |u|^2  |\nabla
u|^2\right) dx.
\end{eqnarray}

Note that there are some terms such as $\int \rho |\dot{u}|^2\, dx $ whose estimates are not clear.  These terms look less frightening than $\int \theta|\nabla u|^2\,  dx$, if one is 
familiar with the techniques for regularity estimates of compressible Navier-Stokes equation. One standard way is to  multiply $(\ref{1.1})_2$ by $u_t$ and integrate the resulting
equation over $\Omega$. Then by the Young's inequality, we obtain that 
\begin{eqnarray}\label{3-3-8}
&&\f12\f{d}{dt}\int\left[ \m|\nabla u|^2+(\m+\l)(\mbox{div}u)^2\right]dx+\int\r|\dot{u}|^2\, dx\nonumber\\
= && \int\r\dot{u}\cdot(u\cdot\nabla){u} \, dx+\int P\mbox{div}u_t \, dx\nonumber\\
\leq &&
\f14\int\r|\dot{u}|^2\, dx+C\int|u|^2|\nabla u|^2\, dx
+\f{d}{dt}\int P\mbox{div}u\, dx-\int
P_t\mbox{div}u\, dx .
\end{eqnarray}

To deal with the term $\int P_t {\rm div} u \,  dx$,  we employ some technique which is a combination of those from \cite{SWZ-JMPA} and \cite{Wen-Zhu}. First, we split $u$ into two parts, $v$ and $w$. 
Let $v = \mathcal{L}^{-1} \nabla P$ and $w= u-v$. (In what follows, $v$ and $w$ always denote $\mathcal{L}^{-1} \nabla P$ and $u-v$. ) As noted in \cite{SWZ-JMPA}, ${\rm div}w$ acts as the effective viscous flux for the bounded domain case. Now 
$$\int P_t {\rm div} u \, dx = \int P_t {\rm div} v \, dx + \int P_t {\rm div} w \, dx .$$
Herein 
\be \label{3-3-8-1} \ba
\int P_t {\rm div} v \, dx  & = - \int \nabla P_t v \, dx \\
& = - \int (\mathcal{L} v)_t v \, dx \\
& = \frac12 \frac{d}{dt} \int  | (-\mathcal{L})^{1/2} v |^2 \, dx,  
\ea \ee
and according to \eqref{3-2-2},
\be \label{3-3-8-2} \ba
& \int P_t  {\rm div} w \, dx  \\= & \frac{R}{c_v} \left[ \int (\rho E)_t {\rm div} w \, dx -  \int \frac12 (\rho |u|^2)_t {\rm div} w \, dx \right] \\
= & \frac{R}{c_v} \left\{ \int \left( \rho E u + P u - \kappa \nabla \theta - \mu \nabla u\cdot u - \mu u\cdot \nabla u  - \lambda u{\rm div}u  \right) \cdot \nabla {\rm div} w \, dx  \right. \\
&  \left.  + \frac12\int  {\rm div} (\rho u) |u|^2 {\rm div} w \, dx        - \int \rho u_t \cdot u  {\rm div} w \, dx      \right\} \\
= & \frac{R}{c_v} \left\{ \int \left( ( c_{\nu}+R) \rho \theta u + \frac12 \rho |u|^2  u  - \kappa \nabla \theta - \mu \nabla u \cdot u - \mu u \cdot \nabla u - \lambda u {\rm div} u    \right) \cdot \nabla {\rm div}w \, dx  \right.\\
& \left.    - \int \frac12  \rho  |u|^2 u\cdot \nabla {\rm div} w \, dx -  \int \rho \dot{u} \cdot u {\rm div} w \, dx \right\}\\
= & \frac{R}{c_v} \left\{\int  \left[ ( c_{v}+R) \rho \theta u - \kappa \nabla \theta - \mu \nabla u\cdot u - \mu u \cdot \nabla u - \lambda u {\rm div} u    \right]\cdot \nabla {\rm div}w \, dx \right.
\\ &\  \left.- \int \rho \dot{u} \cdot u {\rm div} w\, dx  \right\}.
\ea \ee
By virtue of Proposition \ref{prop2.1}, we have 
\be \label{3-3-8-3}
\|\nabla v\|_{L^2} \leq C \| \rho \theta\|_{L^2},
\ee
and 
\be \label{3-3-8-4}
\|\nabla^2 w\|_{L^2} \leq C \|\rho \dot{u}\|_{L^2}. 
\ee
Making use of these two inequalities \eqref{3-3-8-3} and \eqref{3-3-8-4}, we obtain that
\be \label{3-3-8-5} \ba
-\int P_t {\rm div} w \, dx &  \leq C \left(\|\sqrt{\rho} \theta\|_{L^2} \|u\|_{L^\infty} + \|\nabla \theta\|_{L^2} + \|\nabla u\|_{L^2}\|u\|_{L^\infty}\right) \|\nabla {\rm div} w\|_{L^2} \\
&\ \ \ + C \|\sqrt{\rho} \dot{u}\|_{L^2} \|u\|_{L^\infty} \|{\rm div} w\|_{L^2}\\
& \leq C_\delta \|\sqrt{\rho } \theta\|_{L^2}^2 \|u\|_{L^\infty}^2 + C_\delta \|\nabla \theta\|_{L^2}^2 + C_\delta \|\nabla u\|_{L^2}^2 \|u\|_{L^\infty}^2 + \delta \|\sqrt{\rho} \dot{u}\|_{L^2}^2. \\
\ea \ee

Substituting \eqref{3-3-8-1} and \eqref{3-3-8-5}  into \eqref{3-3-8},  we
obtain after choosing $\delta$ suitably small that
\begin{eqnarray} \label{3-3-10}
&&\f{d}{dt}\int\left(\m|\nabla u|^2+(\m+\l)(\mbox{div}u)^2-
2P\mbox{div}u +   |(-\mathcal{L})^{1/2} v|^2 
\right)dx+ \int\r|\dot{u}|^2\, dx\nonumber\\
&& \leq C \|\sqrt{\rho} \theta\|_{L^2}^2 \|u\|_{L^\infty}^2 + C_1 \kappa \|\nabla \theta\|_{L^2}^2 + C \|\nabla u\|_{L^2}^2 \|u\|_{L^\infty}^2~. \end{eqnarray}

Choose some  constant $C_2\ge C_1+1$  suitably large such that
\be\nonumber  \mu|\nabla u|^2-2P{\rm div} u+C_2\rho\t^2\ge
\frac{3 \mu }{4}|\nabla u|^2+\rho\t^2.\ee 
Adding   \eqref{3-3-7}
multiplied by $C_2$ to \eqref{3-3-10}, we have after choosing $\eta$ suitably small that
\begin{eqnarray}\label{3-3-11}
&&\f{d}{dt}\int\left(\m|\nabla u|^2+(\m+\l)(\mbox{div}u)^2-
2P\mbox{div}u +C_2\r\t^2 +  |(-\mathcal{L})^{1/2} v|^2 
 \right) dx  \nonumber\\
&&\quad+ \f12\int\r|\dot{u}|^2\, dx+ \mu \int\t|\nabla
u|^2dx+\k \int|\nabla\t|^2\, dx\nonumber\\
 && \ \ \leq  C \|{\rm div}u\|_{L^\infty} \int \rho \theta^2 \, dx + C \|u\|_{L^\infty}^2 \left(\int \rho \theta^2 \, dx + \|\nabla u\|_{L^2}^2\right)~.
\end{eqnarray}

Let \be \nonumber \Psi(t) = e+ \sup_{\tau\in [0, t]} \left( \|\nabla u(\tau)\|_{L^2}^2 + \int \rho \theta^2(\tau)\,  dx \right) + \int_0^t \int \left( 
\rho|\dot{u}|^2 + \theta|\nabla u|^2 + |\nabla \theta|^2 \right) dx d\tau.\ee
By virtue of Gownwall's inequality, for every $0\leq s \leq T < T^*$, 
\be \label{3-3-12}
\Psi(T) \leq C \Psi(s) \exp \left\{ C \int_s^T \|u(\tau)\|_{L^\infty}^2 \, d\tau \right\}.
\ee

 Now it is time to get a good control of $\|u\|_{L^2(s, T; L^\infty)}$. Making use of Lemma 2.3, we can get that
\be\label{3-3-13} \ba
\|u\|_{L^2(s, T; L^\infty)}^2 &  \leq C \left( 1+ \|u\|_{L^2(s, T; H^1)}^2 \ln \left( e+ \|u\|_{L^2(s, T; W^{1,3}) }\right) \right).
\ea \ee
By Proposition \ref{prop2.1} and Lemma \ref{lem3.4},
\be \nonumber \ba
\|u\|_{W^{1,3}} & \leq C \|w\|_{W^{2, \frac65}} + \|v\|_{W^{1,3}}\\
&\leq C \|\rho \dot{u}\|_{L^{\frac65}} + C \|P \|_{L^3} + C\|u\|_{L^2}\\
& \leq C \| \rho \dot{u}\|_{L^2 } + C \|\nabla \theta \|_{L^2} + C\|\nabla u\|_{L^2} + C,
\ea \ee
which implies that
\be \label{3-3-15} \ba
\ \ \ \|u\|_{L^2(s, T; W^{1,3}) } & \leq C \|\rho \dot{u} \|_{L^2(s, T; L^2)} + C \|\nabla \theta \|_{L^2(s, T; L^2)} +  C \|\nabla u\|_{L^2(s, T; L^2)} + C \\
& \leq C \Psi(T) .
\ea \ee
Substituting \eqref{3-3-15} to \eqref{3-3-13}, 
\be\label{3-3-16}
\|u \|_{L^2(s, T; L^\infty)}^2 \leq  C \left( 1 + \| u \|_{L^2(s,T; H^1)}^2 \ln \left(     C \Psi(T) \right) \right).
\ee

Taking this inequality \eqref{3-3-16} back to \eqref{3-3-12}, then we get 
\be \nonumber
\Psi(T) \leq  C \Psi(s)  (C \Psi(T) )^{ C_3 \|u\|_{L^2(s, T; H^1)}^2 }~.
\ee
Recalling the energy like estimate \eqref{energy-estimate}, we choose some $s$ which is close enough to  $T^*$ such that  
$$\lim_{T\rightarrow T^*-} C_3 \|u\|_{L^2(s, T; H^1)}^2 \leq \frac12,$$  then
\be \label{3-3-19}
\Psi(T) \leq C \Psi(s)^2< \infty,
\ee
which completes the proof for Lemma \ref{lem3.3}.

\vspace{4mm}Lemma \ref{lem3.3} tells that 
\be \nonumber
\lim_{T\rightarrow T^*-} \|\nabla u\|_{L^\infty(0, T; L^2)} <  \infty~,
\ee
which implies that 
\be \nonumber
\lim_{T\rightarrow T^*-} \|u\|_{L^4 (0, T; L^4)} <  \infty~.
\ee
According Huang-Li-Wang \cite{HLW}'s criterion \eqref{HLW}, we can claim here that the strong solution can be extended. For readers' convenience, we give the complete proof. 
The remaining proof consists of higher order estimates of the solutions which are needed to
guarantee the extension of local strong solution to be a global one
under the conditions \eqref{criterion-1} and \eqref{3-3-1}. Compared to \cite{HLW}, there are some slight changes, since we consider the bounded case, instead of the whole space one.

\begin{lemma}\label{lem3.5}
Under the condition \eqref{criterion-1}, it holds that for $0\leq
T<T^\ast$,
\begin{equation}\label{3-5-1}
\sup_{0\leq t\leq
T}\int\left(|\nabla\t|^2+\r|\dot{u}|^2\right)dx+\int_{0}^{T}\int\left(\r\dot{\t}^2+|\nabla\dot{u}|^2\right)dxdt\leq
C.
\end{equation}

\end{lemma}

 {\it Proof}. First, applying $\dot{u}^j[\partial_t+\mbox{div}(u\cdot)]$ to $ (\ref{1.1})_2^j$ and integrating the resulting
equation over $\Omega$, we obtain after integration by parts that
\begin{eqnarray}\label{3-5-2}\f12
\f{d}{dt}\int\r|\dot{u}|^2\, dx&=&-\int\dot{u}^j[\partial_jP_t+\mbox{div}(u\partial_jP)]
\, dx+\m\int\dot{u}^j[\Delta{u}^j_t+\mbox{div}(u\Delta{u}^j)]\, dx\nonumber\\
&&+(\m+\l)\int\dot{u}^j[\partial_j\mbox{div}u_t+\mbox{div}(u\partial_j{\rm div}u)] \, dx\nonumber\\
&=&\sum_{i=1}^{3}N_i.
\end{eqnarray}
We get  after integration by parts and
using  the equation $(\ref{1.1})$ that \be\la{3-5-3} \ba
N_1 & = - \int   \dot{u}^j[\partial_jP_t + \div (\p_jPu)]\, dx \\
& =R\int   \p_j\dot{u}^j\left(\n \dot\te-\r u\cdot\na\te-\te
 u\cdot\na\n-\te  \n\div u\right) dx + \int  \p_k\dot{u}^j\p_jPu^k \, dx  \\
& =R\int   \p_j\dot{u}^j\cdot \n \dot\te \, dx
   -  \int    P\pa_k\dot{u}^j\p_ju^k\,  dx
 \\&\le \frac{ \mu}{8} \|\nabla\dot{u}\|_{L^2}^2
 +C    \|\n \dot\te\|_{L^2}^2
+C\int\n^2\te^2|\na u|^2 \, dx   \\&\le \frac{ \mu}{8} \|\nabla\dot{u}\|_{L^2}^2 + C \|\rho \dot{\theta} \|_{L^2}^2 + 
C \|\theta\|_{L^4}^4 + C \|\nabla u\|_{L^4}^4 .\ea \ee
Integration by parts leads to
\be\label{3-5-4} \ba N_2 & =  \mu\int
 \dot{u}^j[\triangle u_t^j
 + \div (u\triangle u^j)]\, dx \\
& = - \mu\int  \left(\p_i\dot{u}^j\p_iu_t^j +
\triangle u^ju\cdot\nabla\dot{u}^j\right)dx \\
& = -  \mu\int \left(|\nabla\dot{u}|^2 -
\p_i\dot{u}^ju^k\p_k\p_iu^j - \p_i\dot{u}^j\p_iu^k\p_ku^j +
\triangle u^ju\cdot\nabla\dot{u}^j\right) dx \\
& = - \mu\int  \left(|\nabla\dot{u}|^2 + \p_i\dot{u}^j
\p_iu^j\div u - \p_i\dot{u}^j\p_iu^k\p_ku^j - \p_iu^j\p_iu^k\p_k\dot{u}^j
\right) dx \\
&\le -\frac{7\mu}{8} \int |\nabla\dot{u}|^2\, dx  + C \int
 |\nabla u|^4dx  . \ea \ee
Similarly, we have
\begin{eqnarray}\label{3-5-5}
N_3&\leq&
-\f78(\m+\l)\|\mbox{div}\dot{u}\|_{L^2}^2+C\int|\nabla{u}|^4\, dx.
\end{eqnarray}
Substituting \eqref{3-5-3}-\eqref{3-5-5} into \eqref{3-5-2} implies
\begin{equation}\label{3-5-6}\ba
 \f{d}{dt}\int\r|\dot{u}|^2\, dx+\m\|\nabla\dot{u}\|_{L^2}^2 &\leq
C\int\r\dot{\t}^2\, dx+C\|\theta \|_{L^4}^4+C \int|\nabla{u}|^4\, dx 
\\ &\leq
C\int\r\dot{\t}^2\, dx+C\|\nabla\t\|_{L^2}^4+C\|\sqrt{\rho} \dot{u}\|^{4}_{L^2}+ C,\ea
\end{equation}
where for the last inequality we have used the fact, 
\be\label{3-5-6-add} \ba 
\|\nabla u\|_{L^4} \leq \|\nabla v\|_{L^4} + \|\nabla w\|_{L^4} & \leq C \|\rho \theta\|_{L^4} + C \|\rho \dot{u}\|_{L^{4/3}} \\
&\leq C \|\nabla \theta\|_{L^2} + C \|\sqrt{\rho} \dot{u}\|_{L^2} + C~,
\ea \ee
owing to  Proposition 2.1 and  Lemma \ref{lem3.4}.

Next, multiplying $(\ref{1.1})_3 $ by $  \dot\te$ and integrating
the resulting equation over $\Omega $ yield  that
 \be\la{3-5-7} \ba  c_v \int\rho|\dot{\te}|^2dx +\frac{\ka}{2} \frac{d}{dt}  \int |\na\theta|^2\, dx 
&= \ka \int \Delta \theta \cdot (u \cdot \nabla \theta) \, dx
+\lambda \int  (\div u)^2\dot\te  \, dx\\&
\quad+2\mu \int |\mathfrak{D}(u)|^2\dot\te \, dx -R \int\n\te \div
u\dot\te \, dx \\&\triangleq \sum_{i=1}^4 I_i . \ea\ee

We estimate each $I_i (i=1,\cdots,4)$ as follows:

First, it follows from Sobolev embedding theory that for any $\epsilon \in(0,1],$
 \be\la{3-5-8} \ba\int \te^2|\na u|^2dx&\le C\|\te\|_{L^\infty}^2\|\na u\|_{L^2}^2\\ 
&\le \epsilon \|\na^2\te\|_{L^2}^2+C_\epsilon \|\na\te\|_{L^2}^2+ C_\epsilon ,\ea\ee
 which together with the
standard $W^{2,2}$-estimate of $(\ref{1.1})_3$ gives
\bnn \la{3-5-8-1} \ba
\|\t\|_{H^2 }^2 &\leq C\int\r\dot{\t}^2\, dx +C\int\r^2\t^2|\nabla
u|^2\, dx+C\int|\nabla u|^4\, dx + C \|\theta\|_{L^2}^2  \\ 
&\leq C\int\r\dot{\t}^2\, dx +C \epsilon \|\na^2\te\|^2_{L^2} + C_\epsilon \|\na\te\|^2_{L^2} + C\int|\nabla u|^4\, dx + C_\epsilon  .\ea
\enn 
Hence, choosing some $\epsilon$ small enough, we have
\be\label{3-5-9}
\|\t\|_{H^2}^2 \leq C\int\r\dot{\t}^2\, dx  +C\|\na\te\|^2_{L^2}+C\int|\nabla u|^4\, dx + C .
\ee
Consequently, by Gagliardo-Nirenberg inequality, 
\be\la{3-5-10} \ba  |I_1|&\le   C \int |\nabla \theta| |\nabla^2 \theta| |u| dx \\
&\le C \|\nabla^2 \theta\|_{L^2} \|\nabla \theta\|_{L^4} \|u\|_{L^4}  \\
&\le C \|\nabla^2 \theta\|_{L^2} \left(  \|\nabla^2 \theta\|_{L^2} + \|\nabla \theta\|_{L^2}  \right)^{1/2} \|\nabla \theta\|_{L^2}^{1/2}\\
&\le \de  \|\na^2\te\|^2_{L^2} +C_\de
\|\na\te\|^2_{L^2} \\
&\leq C\de\int\r\dot{\t}^2\, dx +C_\de\|\na\te\|^2_{L^2}+C_\de\int|\nabla u|^4\, dx + C_\delta .\ea\ee

Next, integration by parts yields that, for any $\eta,\de\in (0,1],$
\be\la{3-5-11}\ba I_2 =&\lambda \int (\div u)^2 \te_t
\, dx+\lambda \int (\div u)^2u\cdot\na\te
\, dx\\=&\lambda \left(\int (\div u)^2 \te
\, dx\right)_t-2\lambda  \int \te \div u \div (\dot u-u\cdot\na u)
\, dx \\
& \ \ +\lambda \int (\div u)^2u\cdot\na\te
\, dx\\=&\lambda \left(\int (\div u)^2 \te
\, dx\right)_t-2\lambda \int \te \div u \div \dot
u\, dx\\&\ \ +2\lambda  \int \te \div u \div ( u\cdot\na u) \, dx
+\lambda \int (\div u)^2u\cdot\na\te \, dx\\= &\lambda
\left(\int (\div u)^2 \te dx\right)_t-2\lambda \int \te \div u
\div \dot u\, dx 
\\&\ \ +2\lambda \int \te \div u \pa_i u^j\pa_j  u^i \, dx
+ \lambda \int u \cdot\na\left(\te   (\div u)^2 \right) dx
 \\ \le &\lambda\left( \int (\div u)^2 \te \, dx\right)_t  +C\left\|\te|\na u| \right \|_{L^2}\left(\|\na \dot u\|_{L^2}+\|\na u\|_{L^4}^2\right)  \\ 
\le &\lambda\left( \int (\div u)^2 \te \, dx\right)_t  + \eta\|\na \dot u\|_{L^2}^2+ C \de\int\n \dot\te^2 \, dx +C_{\de,\eta}\|\na\te\|_{L^2}^2 +C_{\de}\|\na u\|_{L^4}^4 + C_{\delta, \eta}
,\ea\ee where in the last inequality, we have used \eqref{3-5-8} and \eqref{3-5-9}.

 Then, similar to \eqref{3-5-11}, we have that, for
any $\eta,\de \in (0,1],$
\be \la{3-5-12}\ba I_3\le& 2\mu\left( \int
|\mathfrak{D}(u)|^2 \te \, dx\right)_t+\eta\|\na \dot u\|_{L^2}^2+ C\de\int\n \dot\te^2\, dx \\
&+C_{\de,\eta}\|\na\te\|_{L^2}^2 +C_{\de}\|\na u\|_{L^4}^4+ C_{\delta, \eta}. \ea\ee

Finally, it follows from \eqref{3-5-8} and \eqref{3-5-9} that
 \be\la{3-5-13}\ba
 |I_4|   & \le   \de   \int \n \dot\te^2\, dx +C_\de  \int
\te^2|\na u|^2 \, dx \\ & \le C\de\int\r\dot{\t}^2 \, dx +C_\de\|\na\te\|^2_{L^2} + C_\delta.  \ea\ee

Substituting \eqref{3-5-10}-\eqref{3-5-13} into \eqref{3-5-7}, we obtain after choosing $\de$ suitably small that, for any $\eta\in (0,1],$
\begin{eqnarray}\label{3-5-14}
&&\f{d}{dt}\int\left(\f\k2|\nabla\t|^2-\t\left[\l(\mbox{div}
u)^2+2\m|\mathfrak
{D}(u)|^2\right]\right)dx+\frac{c_v}{2}\int\r\dot{\t}^2 \, dx\nonumber\\
&&\ \ \leq
C\eta\|\nabla\dot{u}\|_{L^2}^2+C_\eta\|\nabla{u}\|_{L^4}^4+C_\eta
\|\nabla\t\|_{L^2}^2 + C_\eta \nonumber\\
&&\ \ \leq
C\eta\|\nabla\dot{u}\|_{L^2}^2+C_\eta\|\sqrt{\rho} \dot u\|_{L^2}^4+C_\eta
\|\nabla\t\|_{L^2}^4+C_\eta~,
\end{eqnarray}
where the last inequality is owing to \eqref{3-5-6-add}.

Noticing that 
\be \ba 
& \int \theta \left[ \lambda ({\rm div} u)^2 + 2\mu |D(u)|^2  \right]dx \\
&\ \ \leq  C \|\theta\|_{L^6} \|\nabla u\|_{L^{12/5}}^2 \\
&\ \ \leq C \left(\|\nabla \theta\|_{L^2} + 1\right) \cdot \|\nabla u\|_{L^2}^{4/3} \|\nabla u\|_{L^4}^{2/3} \\
&\ \ \leq  C  \left(\|\nabla \theta\|_{L^2} + 1\right) \left( \|\sqrt{\rho}\dot{u} \|_{L^2}^{2/3} + \|\rho \theta\|_{L^4}^{2/3}  + 1 \right)\\
&\ \ \leq  \frac{\kappa}{4}\|\nabla \theta\|_{L^2}^2 + \eta^{1/2} \|\sqrt{\rho}\dot{u}\|_{L^2}^2 + C_\eta,
\ea \ee
so adding \eqref{3-5-6} multiplied by $2\eta^{1/2}$ to  \eqref{3-5-14}, we obtain \eqref{3-5-1} after choosing $\eta$ suitably small and  using Gronwall's inequality.
 Thus we complete the proof of Lemma \ref{lem3.5}.

As a corollary, we can bound $\|\theta\|_{L^4}$ and $\|\nabla u\|_{L^4}$. 
\begin{corollary} \label{cor3.5}
Under the condition \eqref{criterion-1}, it holds that for $0\leq
T<T^\ast$,
\begin{equation}\label{3-5-16}
\sup_{0\leq t\leq
T} \left( \|\theta\|_{L^4} + \|\nabla u\|_{L^4} \right) \leq 
C.
\end{equation}
\end{corollary}

{\it Proof}. By virtue of Lemma \ref{lem3.4} and Lemma \ref{lem3.5},
\be \nonumber
\|\theta\|_{L^4}  \leq C \|\nabla \theta\|_{L^2} + C \leq C. \ee

Consequently, according to \eqref{3-5-6-add} and Lemma \ref{lem3.5},
\be \nonumber
\|\nabla u\|_{L^4}  \leq C \|\rho \dot{u} \|_{L^2 } + C\|\rho \theta\|_{L^4} + C \|\nabla u\|_{L^2} + C  \leq  C.
\ee

\vspace{4mm}

Next, we will derive the desired estimates for $\dot\t$. In fact, we
have
\begin{lemma}\label{lem3.6}
Under the condition \eqref{criterion-1}, it holds that for $0\leq
T<T^\ast$,
\begin{equation}\label{3-6-1}
\sup_{0\leq t\leq
T}\int\rho \dot{\t}^2\, dx+\int_{0}^{T}\|\nabla{\dot{\t}}\|_{L^2}^2\, dt\leq
C.
\end{equation}
\end{lemma}

 {\it Proof}. Applying the operator $\pa_t+\div(u\cdot) $ to (\ref{1.1})$_3$, and using (\ref{1.1})$_1$, we get
\be\la{3-6-2}\ba
&c_v \n \left(\pa_t\dot \te+u\cdot\na\dot \te\right)\\
&=\ka \Delta \dot \te+\ka\left( \div u\Delta \te  -
\pa_i\left(\pa_iu\cdot\na \te \right)- \pa_iu\cdot\na \pa_i\te
\right)  -R\n \dot\te \div u-R\n \te\div \dot u\\
&\quad +\left( \lambda (\div u)^2+2\mu |\mathfrak{D}(u)|^2\right)\div u
 +2\lambda \left( \div\dot u-\pa_ku^l\pa_lu^k\right)\div u
\\&\quad
+ \mu (\pa_iu^j+\pa_ju^i)\left( \pa_i\dot u^j+\pa_j\dot
u^i-\pa_iu^k\pa_ku^j -\pa_ju^k\pa_ku^i\right).\ea\ee 

Multiplying (\ref{3-6-2}) by $\dot \te,$  we obtain after integration by parts  and Corollary \ref{cor3.5} that 
\be\la{3-6-3}\ba &
\frac{c_v}{2}\left(\int \n |\dot\te|^2dx\right)_t + \ka
\|\na\dot\te\|_{L^2}^2 \\ \le & 
 C  \int|\na u|\left(|\na^2\te||\dot\te|+ |\na \te| |\na
\dot\te|\right)dx + C   \int|\na
u|^2|\dot\te|  |\na u| \, dx
 \\&\quad +C   \int  \n  |\dot
\te|^2|\na u| \, dx   +C \int  \n \te  |\na\dot u| |\dot
\te| \, dx +C   \int  |\na u| |\na\dot u| |\dot
\te| \, dx \\
\leq  & C \|\nabla u\|_{L^4} \|\nabla^2 \theta \|_{L^2} \|\dot{\theta}\|_{L^4} + 
C \|\nabla u\|_{L^4} \|\nabla \theta\|_{L^2} \|\nabla \dot{\theta}\|_{L^2}  \\ 
& \ \ + C \|\nabla u\|_{L^4}^3 \|\dot{\theta}\|_{L^4} + C \|\nabla u\|_{L^4} \|\sqrt{\rho}\dot{\theta} \|_{L^2} \|\dot{\theta}\|_{L^4} \\
& \ \ + C \|\sqrt{\rho}\dot{\theta}\|_{L^2} \|\nabla \dot{u}\|_{L^2} \|\theta \|_{L^\infty} + C \|\nabla u\|_{L^4} \| \nabla \dot{u}\|_{L^2}
\|\dot{\theta}\|_{L^4}    \\
\leq & \frac{\kappa}{2} \|\nabla \dot{\theta}\|_{L^2}^2 + C \|\nabla^2 \theta\|_{L^2} \|\dot{\theta}\|_{L^4} + C \|\dot{\theta}\|_{L^4} + C \|\sqrt{\rho}\dot{\theta}\|_{L^2} \|\dot{\theta}\|_{L^4} \\
& \ \ + C \|\sqrt{\rho}\dot{\theta}\|_{L^2} \|\nabla \dot{u}\|_{L^2} \|\theta \|_{L^\infty} + C \| \nabla \dot{u}\|_{L^2} \|\dot{\theta}\|_{L^4} + C.   \\  \ea\ee

It follows from \eqref{3-5-9}  and Lemma \ref{lem3.5} that
\be \la{3-6-4} \ba
\|\nabla^2\theta\|_{L^2} & \leq C \|\sqrt{\rho}\dot{\theta}\|_{L^2} + C \|\nabla \theta\|_{L^2} + C \|\nabla u\|_{L^4}^2 + C \\
& \leq C \|\sqrt{\rho} \dot{\theta}\|_{L^2} + C.
\ea \ee
For the estimate for $\|\dot{\theta}\|_{L^4}$, we will follow the method used in Lemma \ref{lem3.4}. Let $\bar{\dot{\theta}} = \frac{1}{|\Omega|}\int \dot{\theta} \, dx $, 
\be\la{3-6-5} \ba 
\bar{\dot{\theta}} \int \rho \, dx & \leq \left|\int \rho \left(\dot{\theta} - \bar{\dot{\theta} } \right) dx  \right| + \left|  \int \rho \dot{\theta} \, dx \right| \\
& \leq C \|\nabla \dot{\theta}\|_{L^2} +  C \|\sqrt{\rho} \dot{\theta }\|_{L^2}, \ea
\ee
which together with Poincar\'e's inequality leads to 
\be \la{3-6-6} \ba 
\|\dot{\theta}\|_{L^4} & \leq  C \|\nabla \dot{\theta}\|_{L^2} + C \left | \bar{\dot{\theta}}  \right| \\
& \leq  C \|\nabla \dot{\theta}\|_{L^2} + C \|\sqrt{\rho}\dot{\theta}\|_{L^2}. \ea
\ee
And $\|\theta\|_{L^\infty}$ can be estimated as follows, 
\be \la{3-6-7} \ba
\|\theta\|_{L^\infty} &  \leq C \|\nabla^2 \theta\|_{L^2} + C \|\theta\|_{L^2} \\
& \leq C \|\sqrt{\rho} \dot{\theta}\|_{L^2} + C.
\ea
\ee

Substituting \eqref{3-6-4}-\eqref{3-6-7} to \eqref{3-6-3}, we arrive at 
\be \la{3-6-8} \ba 
& c_v \left(\int \n |\dot\te|^2 \, dx\right)_t + \ka
\|\na\dot\te\|_{L^2}^2 \\ \le & 
C \int \rho |\dot{\theta}|^2 \, dx + C \|\sqrt{\rho} \dot{\theta}\|_{L^2} \|\nabla \dot{\theta}\|_{L^2} + C\|\nabla \dot{\theta}\|_{L^2} + C \|\sqrt{\rho} \dot{\theta}\|_{L^2} \\
& \ + C \|\nabla \dot{u}\|_{L^2} \int \rho |\dot{\theta}|^2 \, dx +  C \|\nabla \dot{u}\|_{L^2}  \|\sqrt{\rho} \dot{\theta}\|_{L^2}+ C \|\nabla \dot{u}\|_{L^2}^2 + \frac14 \|\nabla \dot{\theta}\|_{L^2}^2+ C,  \ea
\ee
which together with the Gronwall's inequality completes the proof for Lemma \ref{lem3.6}.

\vspace{3mm}As a corollary, the bounds for $\|\theta\|_{H^2}$ and $\|\theta\|_{L^\infty}$  can be derived.
\begin{corollary}\label{cor3.6-1}
Under the condition \eqref{criterion-1}, it holds that for $0\leq T < T^*$, 
\be \la{3-6-9}
\sup_{0\leq t \leq T} \left( \|\theta\|_{H^2} + \|\theta\|_{L^\infty} \right) \leq C.
\ee
\end{corollary}

{\it Proof}.\ First, it follows from \eqref{3-5-9}, Lemma \ref{lem3.5}, Corollary \ref{cor3.5} and Lemma \ref{lem3.6}
that 
\be \la{3-6-11}
\|\nabla^2 \theta\|_{L^2} \leq C.
\ee
Hence, 
\be \la{3-6-12}
\|\theta\|_{L^\infty} \leq C\|\theta\|_{H^2} \leq C.
\ee

\vspace{3mm}Up to now, we have get the bounds for $\|\rho\|_{L^\infty}$ and $\|\theta\|_{L^\infty}$, which imply other necessary high order estimates for the extension of the strong solution, according to the theorem proved in \cite{FZZ}. We sketch the proof for completeness. 

\vspace{3mm}\begin{corollary}\label{cor3.6-2}Under the condition \eqref{criterion-1}, it holds that for $0\leq T<T^*$,
\be \label{3-6-13}
\sup_{0\leq t \leq T} \|w\|_{H^2} + \int_0^T \left(\|\nabla^2 w\|_{L^p}^2 + \|\nabla w\|_{L^\infty}^2 \right) dt \leq C, \ \ \ p\in (2, \infty). 
\ee
\end{corollary}

{\it Proof}.\ By virtue of Proposition \ref{prop2.1} and Lemma \ref{lem3.5}, 
$$\|w\|_{H^2} \leq C \|\rho \dot{u}\|_{L^2} \leq C, $$
and by Sobolev embedding inequality,
$$\ba 
\|\nabla w\|_{L^\infty} \leq  C\|\nabla  w\|_{W^{1, p}}  \leq C \|\rho \dot{u}\|_{L^p}  \leq C \| \nabla \dot{u}\|_{L^2},
\ea $$
which implies \eqref{3-6-13}.


\vspace{3mm} The next lemma is used to bound the density gradient and $\|u\|_{H^2}$.
\begin{lemma}\label{lem3.7}
Under the condition \eqref{criterion-1},  it holds that for $0\leq T < T^*$, 
\begin{equation}\label{3-7-1}
\sup_{0\leq t\leq T}(\|\rho \|_{W^{1, q}}+\|u \|_{H^2})\leq C.
\end{equation}
\end{lemma}

{\it Proof}.  \ For $2\leq p \leq q $, 
$|\nabla\r|^p$ satisfies the following equation
\bnn\label{3-7-2}
&&(|\nabla\r|^p)_t+\mbox{div}(|\nabla\r|^p u)+(p-1)|\nabla\r|^p\mbox{div}u\nonumber\\
&&\qquad
+p|\nabla\r|^{p-2}(\nabla\r)^{tr}\nabla{u}(\nabla\r)+p\r|\nabla\r|^{p-2}
\nabla\r\cdot\nabla\mbox{div}u=0.
\enn
Hence, 
 \be \label{3-7-3} \ba
\f{d}{dt}\|\nabla\r\|_{L^p}  & \leq C(1+\|\nabla
u\|_{L^\infty})\|\nabla\r\|_{L^p}+C\|\nabla^2 u\|_{L^p} \\
& \leq C\left(1+ \|\nabla v\|_{L^\infty} + \|\nabla w\|_{\infty} \right) \|\nabla \rho\|_{L^p} + C \left( \|\nabla^2 v\|_{L^p} + \|\nabla^2 w\|_{L^p}\right)\\
& \leq C \left(1+ \|\nabla v\|_{L^\infty} + \|\nabla w\|_{\infty} \right) \|\nabla \rho\|_{L^p} + C  \|\nabla^2 w\|_{L^p} + C,
\ea\ee
  where for the last inequality we used the fact \be\ba
 \la{3-7-4} \|\nabla^2 v\|_{L^p} & \leq C \|\nabla(\rho \theta) \|_{L^p} \\
& \leq C \|\nabla \rho\|_{L^p} \|\theta\|_{L^\infty} + C \|\nabla \theta\|_{L^p } \|\rho\|_{L^\infty} \\
&\leq C \|\nabla \rho\|_{L^p}  + C .\ea\ee

To bound $\|\nabla v\|_{L^\infty}$, we make use of the endpoint case of Proposition \ref{prop2.1}, Lemma \ref{lem2.2} and \eqref{3-7-4},
\begin{equation}\label{3-7-5} \ba
\|\nabla v\|_{L^\infty} & \leq C \left(1 + \|\nabla v\|_{BMO(\Omega)} \ln (e + \|\nabla v\|_{W^{1, p}}) \right) \\
& \leq C \left(1 + (\|P\|_{L^\infty} + \|P\|_{L^2} ) \ln (e + \|\nabla v\|_{W^{1,p}} ) \right)\\
& \leq C \left(  1 + \ln (e + \|\nabla \rho \|_{L^p} )       \right).
\ea \end{equation}
Substituting \eqref{3-7-5} into \eqref{3-7-3}, we get that
\be \label{3-7-6}\ba
& \frac{d}{dt} \left( e + \|\nabla \rho\|_{L^p} \right)\\
 \leq & C \left( 1 + \|\nabla w\|_{L^\infty} \right) \|\nabla \rho\|_{L^p} + C \ln (e + \|\nabla \rho\|_{L^p} ) \|\nabla \rho\|_{L^p} + C \|\nabla^2 w\|_{L^p}+C,
\ea \ee
which together with Gronwall's inequality and Corollary \ref{cor3.6-2} gives that
\be\la{3-7-7}
\sup_{0<T< T^*} \|\nabla \rho\|_{L^p} \leq C.
\ee
Let $p=q$, then we get the bound of $\|\rho\|_{W^{1,q}}$.

Moreover, Let $p=2$ in \eqref{3-7-7}, then by Corollary \ref{cor3.6-2} and \eqref{3-7-4},
\be\la{3-7-8}
\|\nabla^2 u\|_{L^2} \leq \|\nabla^2 w\|_{L^2} + \|\nabla^2 v \|_{L^2} \leq C + C \|\nabla \rho\|_{L^2} \leq C. 
\ee
It completes the proof of Lemma \ref{lem3.7}.


\vspace{4mm}  
In view of Lemmas \ref{lem3.1}-\ref{lem3.7}, it is enough to
extend the strong solution $(\r,u,\t)$ beyond $t\geq T^\ast$. In
fact, note that the generic constants $C$ in Lemmas
\ref{lem3.1}-\ref{lem3.7} remains uniformly bounded for all
$T<T^\ast$, so the functions
$(\r,u,\t)(x,T^\ast)\triangleq\lim_{t\rightarrow
T^\ast-}(\r,u,\t)(x,t)$ satisfy the conditions imposed on the initial
data at the time $t=T^\ast$. Furthermore, standard
arguments yield that $\r\dot{u},\r\dot{\t}\in C([0,T];L^2)$, which
implies 
\be \nonumber  (\r\dot{u},\r\dot{\t})(x,T^\ast)=\lim_{t\rightarrow
T^\ast-}(\r\dot{u},\r\dot{\t})(x, t)\in L^2. \ee
 Hence, \bnn
&&\m\Delta{u} + (\m+\l)\nabla\mbox{div}u- R\nabla(\r\t)|_{t=T^\ast}=\sqrt{\r}(x,T^\ast)g_1(x),\nonumber\\
&&\k\Delta\t+\f\m2|\nabla{u}+(\nabla{u})^{tr}|^2+\l(\mbox{div}u)^2|_{t=T^\ast}=\sqrt{\r}(x,T^\ast)g_2(x)\nonumber,
\enn 
with \bnn  g_1(x)\triangleq
\begin{cases}
\r^{-1/2}(x,T^\ast)(\r\dot{u})(x,T^\ast),~
\mbox{for}~~x\in\{x|\ \r(x,T^\ast)>0\},\\
0,~~~~~~~~~~~~~~~~~~~~~~~~~~~~\mbox{for}~~x\in\{x|\ \r(x,T^\ast)=0\},
\end{cases}
\enn and \bnn g_2(x)\triangleq
\begin{cases}
\r^{-1/2}(x,T^\ast)[\r\dot{\t}+R\r\t\mbox{div}u](x,T^\ast),~
\mbox{for}~~x\in\{x|\ \r(x,T^\ast)>0\},\\
0,~~~~~~~~~~~~~~~~~~~~~~~~~~~~~~~~~~~~~~~~~~\mbox{for}~~x\in\{x|\ \r(x,T^\ast)=0\},
\end{cases}
\enn satisfying $g_1,g_2\in L^2$ due to Lemma \ref{lem3.7} and Corollary \ref{cor3.6-1}. Thus $(\r,u,\t)(x,T^\ast)$ satisfies \eqref{compatible-1} and
\eqref{compatible-2} also. Therefore, one can take $(\r,u,\t)(x,T^\ast)$ as
the initial data and apply Theorem \ref{local-existence} to extend the local
strong solution beyond $T^\ast$. This contradicts the assumption on
$T^{\ast}$. We thus finish the proof of Theorem \ref{thm1.1}.

\vspace{4mm}
{\bf Acknowledgements} \ The author would like to thank Professor Jing Li and Dr. Xiangdi Huang for introducing this topic and all the helpful discussions. The author is partially supported 
by NSF of China under Grant 11241004.

\end{document}